\documentclass[11pt,a4paper]{amsart} 
\usepackage{amssymb}
\usepackage[all]{xy}
\SelectTips{cm}{}
\calclayout 

\title{An axiomatic characterization of the Gabriel-Roiter measure
} 
\author[Henning Krause]{Henning Krause}
\address{Henning Krause\\ Institut f\"ur Mathematik\\
Universit\"at Paderborn\\ 33095 Paderborn\\ Germany.}
\email{hkrause@math.upb.de}
\subjclass[2000]{Primary: 18E10; Secondary: 06A07, 16G10}

\newtheorem{lem}{Lemma}[section]
\newtheorem{prop}[lem]{Proposition}
\newtheorem{cor}[lem]{Corollary}
\newtheorem{thm}[lem]{Theorem}

\theoremstyle{remark}

\theoremstyle{definition}
\newtheorem{exm}[lem]{Example}

\numberwithin{equation}{section}
\newtheorem{rem}[lem]{Remark}

\newcommand{\smatrix}[1]{\left[\begin{smallmatrix}#1\end{smallmatrix}\right]}

\renewcommand{\leq}{\leqslant}
\renewcommand{\geq}{\geqslant}
\newcommand{\Ch}{\operatorname{Ch}\nolimits}
\newcommand{\ind}{\operatorname{ind}\nolimits}

\newcommand{\card}{\operatorname{card}\nolimits}

\newcommand{\comp}{\mathop{\raisebox{+.3ex}{\hbox{$\scriptstyle\circ$}}}}
\newcommand{\lto}{\longrightarrow}

\def\a{\alpha}
\def\b{\beta}

\def\g{\gamma}
\def\p{\phi}

\def\m{\mu}

\def\la{\lambda}

\def\La{\Lambda}

\def\A{{\mathcal A}}

\def\bbN{\mathbb N}

\def\bbQ{\mathbb Q}
\def\bbR{\mathbb R}

\def\bbZ{\mathbb Z}

\begin{document}

\begin{abstract}
Given an abelian length category $\A$, the Gabriel-Roiter measure with
respect to a length function $\ell$ is characterized as a universal
morphism $\ind\A\to P$ of partially ordered sets. The map is defined
on the isomorphism classes of indecomposable objects of $\A$ and is a
suitable refinement of the length function $\ell$.
\end{abstract}

\maketitle 

In his proof of the first Brauer-Thrall conjecture \cite{Ro}, Roiter
used an induction scheme which Gabriel formalized in his report on
abelian length categories \cite{G}. The first Brauer-Thrall conjecture
asserts that every finite dimensional algebra of bounded
representation type is of finite representation type. Ringel
noticed (see the footnote on p.\ 91 of \cite{G}) that the
formalism of Gabriel and Roiter works equally well for studying the
representations of algebras having unbounded representation type. We
refer to recent work \cite{R1,R2,R3} for some beautiful applications.

In this note we present an axiomatic characterization of the
Gabriel-Roiter measure which reveals its combinatorial nature.  Given
a finite dimensional algebra $\La$, the Gabriel-Roiter measure is
characterized as a universal morphism $\ind\La\to P$ of partially
ordered sets. The map is defined on the isomorphism classes of
finite dimensional indecomposable $\La$-modules and is a suitable
refinement of the length function $\ind\La\to\bbN$ which sends a
module to its composition length.

The first part of this paper is purely combinatorial and might be of
independent interest. We study length functions $\la\colon S\to T$ on
a fixed partially ordered set $S$.  Such a length function takes its
values in another partially ordered set $T$, for example $T=\bbN$. We
denote by $\Ch(T)$ the set of finite chains in $T$, together with the
lexicographic ordering. The map $\la$ induces a new length function
$\la^*\colon S\to \Ch(T)$, which we call chain length function because
each value $\la^*(x)$ measures the lengths $\la(x_i)$ of the elements
$x_i$ occuring in some finite chain $x_1<x_2<\ldots <x_n=x$ of $x$ in
$S$. We think of $\la^*$ as a specific refinement of $\la$ and provide
an axiomatic characterization. It is interesting to observe that this
construction can be iterated.  Thus we may consider $(\la^*)^*$,
$((\la^*)^*)^*$, and so on.

The second part of the paper discusses the Gabriel-Roiter measure for
a fixed abelian length category $\A$, for example the category of
finite dimensional $\La$-modules over some algebra $\La$. For each
length function $\ell$ on $\A$, we consider its restriction to the
partially ordered set $\ind\A$ of isomorphism classes of
indecomposable objects of $\A$.  Then the Gabriel-Roiter measure with
respect to $\ell$ is by definition the corresponding chain length
function $\ell^*$. In particular, we obtain an axiomatic
characterization of $\ell^*$ and use it to reprove Gabriel's main
property of the Gabriel-Roiter measure. Note that we work with a
slight generalization of Gabriel's original definition. This enables
us to characterize the injective objects of $\A$ as those objects
where $\ell^*$ takes maximal values for some length function $\ell$.
This is a remarkable fact because the Gabriel-Roiter measure is a
combinatorial invariant, depending only on the poset of indecomposable
objects and some length function, whereas the notion of injectivity
involves all morphisms of the category $\A$.

\section{Chains and length functions}

\subsection*{The lexicographic order on finite chains}
Let $(S,\leq)$ be a partially ordered set.  A subset $X\subseteq S$ is a {\em chain}
if $x_1\leq x_2$ or $x_2\leq x_1$ for each pair $x_1,x_2\in X$.  For a
finite chain $X$, we denote by $\min X$ its minimal and by $\max X$
its maximal element, using the convention
$$\max \emptyset < x< \min \emptyset\quad\text{for all}\quad x\in S.$$ We
write $\Ch(S)$ for the set of all finite chains in $S$  and let
$$\Ch(S,x):=\{X\in\Ch(S)\mid\max X=x\}\quad\text{for}\quad x\in S.$$
On $\Ch(S)$ we consider the {\em lexicographic order} which is defined
by
$$X\leq Y \quad :\Longleftrightarrow\quad \min(Y\setminus X)\leq
\min(X\setminus Y)\quad\text{for}\quad  X,Y\in\Ch(S).$$ 

\begin{rem}
(1) $X\subseteq Y$ implies $X\leq Y$ for $X,Y\in\Ch(S)$.

(2) Suppose that $S$ is totally ordered. Then $\Ch(S)$ is totally
ordered.  We may think of $X\in\Ch(S)\subseteq\{0,1\}^S$ as a string
of $0$s and $1$s which is indexed by the elements in $S$. The usual
lexicographic order on such strings coincides with the lexicographic
order on $\Ch(S)$.
\end{rem}

\begin{exm}
Let $\bbN=\{1,2,3,\cdots\}$ and $\bbQ$ be the set of rational numbers
together with the natural ordering. Then the map
$$\Ch(\bbN)\lto \bbQ ,\quad X\mapsto \sum_{x\in X}2^{-x}$$ is
injective and order preserving, taking values in the interval
$[0,1]$. For instance, the subsets of $\{1,2,3\}$ are ordered as follows:
$$\{\}<\{3\}<\{2\}<\{2,3\}<\{1\}<\{1,3\}<\{1,2\}<\{1,2,3\}.$$
\end{exm}

We need the following properties of the lexicographic order.

\begin{lem}
\label{le:lex}
Let $X,Y\in\Ch(S)$ and $X^*:=X\setminus\{\max X\}$.
\begin{enumerate}
\item $X^*=\max\{X'\in\Ch(S)\mid X'<X\text{ and }\max X'<\max X\}$.
\item If $X^*< Y$ and $\max X\geq\max Y$, then $X\leq Y$.
\end{enumerate}
\end{lem}
\begin{proof}
(1) Let $X'<X$ and $\max X'<\max X$. We show that $X'\leq X^*$. This
    is clear if $X'\subseteq X^*$. Otherwise, we have
$$\min (X^*\setminus X')=\min (X\setminus X')<\min (X'\setminus X)=\min
(X'\setminus X^*),$$ and therefore $X'\leq X^*$.

(2) The assumption $X^*< Y$ implies by definition
$$\min (Y\setminus X^*)<\min (X^*\setminus Y).$$ We consider two
cases. Suppose first that $X^*\subseteq Y$. If $X\subseteq Y$, then
$X\leq Y$. Otherwise, $$\min (Y\setminus X)<\max X=\min (X\setminus Y)$$
and therefore $X<Y$. Now suppose that $X^*\not\subseteq Y$. We use
again that $\max X\geq\max Y$, exclude the case $Y\subseteq X$, and
obtain
$$\min (Y\setminus X)=\min (Y\setminus X^*)<\min (X^*\setminus Y)=\min
(X\setminus Y).$$ Thus $X\leq Y$ and the proof is complete.
\end{proof}

\subsection*{Length functions}\label{se:lfct}
Let $(S,\leq)$ be a partially ordered set. A {\em length function} on
$S$ is by definition a map $\la\colon S\to T$ into a partially ordered
set $T$ satisfying for all $x,y\in S$ the following:
\begin{enumerate}
\item[(L1)] $x< y$ implies $\la(x)< \la(y)$.
\item[(L2)] $\la(x)\leq\la(y)$ or  $\la(y)\leq\la(x)$.
\item[(L3)] $\la_0(x):=\card\{\la(x')\mid x'\in S\text{ and }x'\leq x\}$
is finite.
\end{enumerate}
Two length functions $\la$ and $\la'$ on $S$ are {\em equivalent} if
$$\la(x)\leq\la(y)\quad\Longleftrightarrow\quad
\la'(x)\leq\la'(y)\quad\text{for all}\quad x,y\in S.$$
Observe that (L2) and (L3) are automatically satisfied if $T=\bbN$.  A
length function $\la\colon S\to T$ induces for each $x\in S$ a map
$$\Ch(S,x)\lto\Ch(T,\la(x)),\quad X\mapsto \la(X),$$ and therefore
the following {\em chain length function}
$$S\lto \Ch(T),\quad x\mapsto \la^*(x):=\max\{\la(X)\mid X\in\Ch(S,x)
\}.$$ Note that equivalent length functions induce equivalent chain
length functions.

\begin{exm}
(1) Let $S$ be a poset such that for each $x\in S$ there is a bound
$n_x\in\bbN$ with $\card X\leq n_x$ for all $X\in\Ch(S,x)$. Then the
map $S\to\bbN$ sending $x$ to $\max\{\card X\mid X\in\Ch(S,x)\}$ is a
length function.

(2) Let $S$ be a poset such that $\{x'\in S\mid x'\leq x\}$ is a finite
chain for each $x\in S$. Then the map $\la\colon S\to\bbN$ sending $x$
to $\card\{x'\in S\mid x'\leq x\}$ is a length function. Moreover,
$\la^*$ is a length function and equivalent to $\la$.

(3) Let $\la\colon S\to \bbZ$ be a length function which satisfies in
addition the following properties of a {\em rank function}:
$\la(x)=\la(y)$ for each pair $x,y$ of minimal elements of $S$, and
$\la(x)=\la(y)-1$ whenever $x$ is an immediate predecessor of $y$ in
$S$. Then $\la^*$ is a length function and equivalent to $\la$.
\end{exm}

\subsection*{Basic properties}
Let $\la\colon S\to T$ be a length function and $\la^*\colon S\to
\Ch(T)$ the induced chain length function. We collect the basic
properties of $\la^*$.

\begin{prop}
\label{pr:def}
Let $x,y\in S$. 
\begin{enumerate}
\item[(C0)] $\la^*(x)=\max_{x'<x}\la^*(x')\cup\{\la(x)\}$.
\item[(C1)] $x\leq y$ implies $\la^*(x)\leq\la^*(y)$.
\item[(C2)] $\la^*(x)=\la^*(y)$ implies $\la(x)=\la(y)$.
\item[(C3)] $\la^*(x')<\la^*(y)$ for all $x'<x$ and $\la(x)\geq\la(y)$
imply $\la^*(x)\leq\la^*(y)$.
\end{enumerate}
\end{prop}

The first property shows that the function $\la^*\colon S\to\Ch(T)$
can be defined by induction on the length $\la_0(x)$ of the elements
$x\in S$.  The subsequent properties suggest to think of $\la^*$ as a
refinement of $\la$.

\begin{proof}
To prove (C0), let $X=\la^*(x)$ and note that $\max X=\la(x)$. The
assertion follows from Lemma~\ref{le:lex} because we have
$$X\setminus\{\max X\}=\max\{X'\in\Ch(T)\mid X'<X\text{ and }\max
X'<\max X\}.$$ Now suppose $x\leq y$ and let $X\in\Ch(S,x)$. Then
$Y=X\cup\{y\}\in\Ch(S,y)$ and we have $\la(X)\leq\la(Y)$ since
$\la(X)\subseteq\la(Y)$. Thus $\la^*(x)\leq \la^*(y)$.  If $\la^*(x)=
\la^*(y)$, then
    $$\la(x)=\max\la^*(x)=\max\la^*(y)=\la(y).$$ To prove (C3), we use
    (C0) and apply Lemma~\ref{le:lex} with $X=\la^*(x)$ and
    $Y=\la^*(y)$. In fact, $\la^*(x')<\la^*(y)$ for all $x'<x$ implies
    $X^*<Y$, and $\la(x)\geq\la(y)$ implies $\max X\geq \max Y$. Thus
    $X\leq Y$.
\end{proof}

\begin{cor}\label{co:len}
Let $\la\colon S\to T$ be a length function. Then the induced map
$\la^*$ is a length function.
\end{cor}
\begin{proof}
(L1) follows from (C1) and (C2). (L2) and (L3) follow from the
corresponding conditions on $\la$.
\end{proof}

\subsection*{An axiomatic characterization}
Let $\la\colon S\to T$ be a length function.  We present an axiomatic
characterization of the induced chain length function $\la^*$.  Thus we
can replace the original definition in terms of chains by three simple
conditions which express the fact that $\la^*$ refines $\la$.

\begin{thm}
\label{se:axiom}
Let $\la\colon S\to T$ be a length function.  Then there exists a map
$\m\colon S\to U$ into a partially ordered set $U$ satisfying for all
$x,y\in S$ the following:
\begin{enumerate}
\item[(M1)] $x\leq y$ implies $\m(x)\leq\m(y)$.
\item[(M2)] $\m(x)=\m(y)$ implies $\la(x)=\la(y)$.
\item[(M3)] $\m(x')<\m(y)$ for all $x'<x$ and $\la(x)\geq\la(y)$ imply
$\m(x)\leq\m(y)$.
\end{enumerate}
Moreover, for any map $\m'\colon S\to U'$ into a partially ordered set
$U'$ satisfying the above conditions, we have
$$\m'(x)\leq\m'(y)\quad\Longleftrightarrow\quad
\m(x)\leq\m(y)\quad\text{for all}\quad x,y\in S.$$
\end{thm}

\begin{proof}
We have seen in Proposition~\ref{pr:def} that $\la^*$ satisfies (M1) -- (M3).
So it remains to show that for any map $\m\colon S\to U$ into a
partially ordered set $U$, the conditions (M1) -- (M3) uniquely
determine the relation $\m(x)\leq\m(y)$ for any pair $x,y\in S$.  We
proceed by induction on the length $\la_0(x)$ of the elements
$x\in S$ and show in each step the following for $S_n=\{x\in
S\mid\la_0(x)\leq n\}$.
\begin{enumerate}
\item[(i)] $\{\m(x')\mid x'\in S_n\text{ and }x'\leq x\}$ is a finite set for all
$x\in S$.
\item[(ii)] (M1) -- (M3) determine the relation $\m(x)\leq\m(y)$ for all
$x,y\in S_n$.
\item[(iii)] $\m(x)\leq\m(y)$ or $\m(y)\leq\m(x)$ for all
$x,y\in S_n$.
\end{enumerate}
For $n=1$ the assertion is clear. In fact, $S_1$ is the set of minimal
elements in $S$, and $\la(x)\geq\la(y)$ implies $\m(x)\leq\m(y)$ for
$x,y\in S_1$, by (M3).  Now let $n>1$ and assume the assertion is true
for $S_{n-1}$.  To show (i), fix $x\in S$. The map
$$\{\m(x')\mid x'\in S_n\text{ and }x'\leq x\}\lto \{\m(x')\mid x'\in
S_{n-1}\text{ and }x'\leq x\}\times\{\la(x')\mid x'\leq x\}$$ sending
$\m(x')$ to the pair $(\max_{y<x'}\m(y),\la(x'))$ is well-defined by
(i) and (iii); it is injective by (M3). Thus $\{\m(x')\mid
x'\in S_n\text{ and }x'\leq x\}$ is a finite set. In order to verify
(ii) and (iii), we choose for each $x\in S_n$ a {\em Gabriel-Roiter
filtration}, that is, a sequence
$$x_1<x_2<\ldots <x_{\g(x)-1}<x_{\g(x)}=x$$ in $S$ such that $x_1$ is
minimal and $\max_{x'<x_{i}}\m(x')=\m(x_{i-1})$ for all $1< i\leq
\g(x)$. Such a filtration exists because the elements $\m(x')$ with
$x'<x$ form a finite chain, by (i) and (iii). Now fix $x,y\in S_n$ and
let $I=\{i\geq 1\mid \m(x_i)=\m(y_i)\}$. We consider $r=\max I$ and
put $r=0$ if $I=\emptyset$. There are two possible cases.  Suppose
first that $r=\g(x)$ or $r=\g(y)$. If $r=\g(x)$, then
$\m(x)=\m(x_r)=\m(y_r)\leq\m(y)$ by (M1).  Now suppose $\g(x)\neq
r\neq\g(y)$. Then we have $\la(x_{r+1})\neq\la(y_{r+1})$ by (M2) and
(M3). If $\la(x_{r+1})>\la(y_{r+1})$, then we obtain
$\m(x_{r+1})<\m(y_{r+1})$, again using (M2) and (M3). Iterating this
argument, we get $\m(x)=\m(x_{\g(x)})<\m(y_{r+1})$. From (M1) we get
$\m(x)<\m(y_{r+1})\leq\m(y)$. Thus $\m(x)\leq\m(y)$ or
$\m(x)\geq\m(y)$ and the proof is complete.
\end{proof}

\begin{cor}
Let $\la\colon S\to T$ be a length function and let $\m\colon S\to U$
be a map into a partially ordered set $U$ satisfying {\rm (M1) --
(M3)}. Then $\m$ is a length function.  Moreover, we have for all
$x,y\in S$
$$\m(x)= \m(y)\quad \Longleftrightarrow\quad
\max_{x'<x}\m(x')=\max_{y'<y}\m(y')\text{ and }\la(x)=\la(y).$$
\end{cor}

\subsection*{Iterated length functions}
Let $\la$ be a length function. Then $\la^*$ is again a length
function by Corollary~\ref{co:len}. Thus we may define inductively
$\la^{(0)}=\la$ and $\la^{(n)}=(\la^{(n-1)})^*$ for $n\geq 1$.  In
many examples, we have that $\la^{(1)}$ and $\la^{(3)}$ are
equivalent. However, this is not a general fact. The author is
grateful to Osamu Iyama for suggesting the following example.

\begin{exm}
The following length functions $\la^{(1)}$ and $\la^{(3)}$ are not
equivalent.
$$\xymatrix@=1.0em{
\la^{(0)}:&4\ar@{-}[d]\ar@{-}[dr]&5\ar@{-}[d]\ar@{-}[dr]&6\ar@{-}[d]
&\la^{(1)}:&3\ar@{-}[d]\ar@{-}[dr]&6\ar@{-}[d]\ar@{-}[dr]&5\ar@{-}[d]
&\la^{(2)}:&6\ar@{-}[d]\ar@{-}[dr]&4\ar@{-}[d]\ar@{-}[dr]&2\ar@{-}[d]\\
&3&2&1&&1&2&4&&5&3&1\\
\la^{(3)}:&3\ar@{-}[d]\ar@{-}[dr]&5\ar@{-}[d]\ar@{-}[dr]&6\ar@{-}[d]
&\la^{(4)}:&6\ar@{-}[d]\ar@{-}[dr]&4\ar@{-}[d]\ar@{-}[dr]&2\ar@{-}[d]\\
&1&2&4&&5&3&1 }
$$
\end{exm}

\section{Abelian length categories}

In this section we recall the definition and some basic facts about
abelian length categories. We fix an abelian category $\A$.

\subsection*{Subobjects}
We say that two monomorphisms $\p_1\colon X_1\to X$ and $\p_2\colon
X_2\to X$ in $\A$ are {\em equivalent}, if there exists an isomorphism
$\a\colon X_1\to X_2$ such that $\p_1=\p_2\comp\a$.  An equivalence
class of monomorphisms into $X$ is called a {\em subobject} of
$X$. Given subobjects $\p_1\colon X_1\to X$ and $\p_2\colon X_2\to X$
of $X$, we write $X_1\subseteq X_2$ if there is a morphism $\a\colon
X_1\to X_2$ such that $\p_2=\p_1\comp\a$.  An object $X\neq 0$ is {\em
simple} if $X'\subseteq X$ implies $X'=0$ or $X'=X$.

\subsection*{Length categories}
An object $X$ of $\A$ has {\em finite length} if it has a finite
composition series $$0=X_0\subseteq X_1\subseteq \ldots \subseteq
X_{n-1}\subseteq X_n=X,$$ that is, each $X_i/X_{i-1}$ is simple.  Note
that $X$ has finite length if and only if $X$ is both artinian (i.e.\
it satisfies the descending chain condition on subobjects) and noetherian
(i.e.\ it satisfies the ascending chain condition on subobjects).  An
abelian category is called a {\em length category} if all objects have
finite length and if the isomorphism classes of objects form a set.

Recall that an object $X\neq 0$ is {\em indecomposable} if
$X=X_1\oplus X_2$ implies $X_1=0$ or $X_2=0$. A finite length object
admits a finite direct sum decomposition into indecomposable objects
having local endomorphism rings. Moreover, such a decomposition is
unique up to an isomorphism by the Krull-Remak-Schmidt Theorem.

\begin{exm}
(1) The finitely generated modules over an artinian ring form a length
category.

(2) Let $k$ be a field and $Q$ be any quiver. Then the finite
    dimensional $k$-linear representations of $Q$ form a length
    category.
\end{exm}

\section{The Gabriel-Roiter measure}

Let $\A$ be an abelian length category. The definition of the
Gabriel-Roiter measure of $\A$ is due to Gabriel \cite{G} and was
inspired by the work of Roiter \cite{Ro}.  We present a definition
which is a slight generalization of Gabriel's original definition.
Then we discuss some specific properties.

\subsection*{Length functions}
A {\em length function} on $\A$ is by definition a map $\ell$ which
sends each object $X\in\A$ to some real number $\ell(X)\geq 0$ such that
\begin{enumerate}
\item[($\ell$1)] $\ell(X)=0$ if and only if $X=0$, and
\item[($\ell$2)] $\ell(X)=\ell(X')+\ell(X'')$ for every exact sequence $0\to
X'\to X\to X''\to 0$.
\end{enumerate}
Note that such a length function is determined by the set of values
$\ell(S)>0$, where $S$ runs through the isomorphism classes of simple
objects of $\A$.  This follows from the Jordan-H\"older Theorem. We
write $\ell_1$ for the length function satisfying $\ell_1(S)=1$ for
every simple object $S$. Observe that $\ell_1(X)$ is the usual
composition length of an object $X\in\A$.

\subsection*{The Gabriel-Roiter measure}
We consider the set $\ind\A$ of isomorphism classes of indecomposable
objects of $\A$ which is partially ordered via the subobject relation
$X\subseteq Y$.  Now fix a length function $\ell$ on $\A$. The map
$\ell$ induces a length function $\ind\A\to\bbR$ satisfying (L1) --
(L3), and the induced chain length function
$\ell^*\colon\ind\A\to\Ch(\bbR)$ is by definition the {\em
Gabriel-Roiter measure} of $\A$ with respect to $\ell$.  Gabriel's
original definition \cite{G} is based on the length function
$\ell_1$. Whenever it is convenient, we substitute $\m=\ell^*$.

\subsection*{An axiomatic characterization}\label{se:GRaxiom}
The following axiomatic characterization of the Gabriel-Roiter measure
is the main result of this note.

\begin{thm}
Let $\A$ be an abelian length category and $\ell$ a length function on
$\A$.  Then there exists a map $\m\colon \ind\A\to P$ into a partially
ordered set $P$ satisfying for all $X,Y\in\ind\A$ the following:
\begin{enumerate}
\item[(GR1)] $X\subseteq Y$ implies $\m(X)\leq\m(Y)$.
\item[(GR2)] $\m(X)=\m(Y)$ implies $\ell(X)=\ell(Y)$.
\item[(GR3)] $\m(X')<\m(Y)$ for all $X'\subset X$ and $\ell(X)\geq\ell(Y)$ imply
$\m(X)\leq\m(Y)$.
\end{enumerate}
Moreover, for any map $\m'\colon \ind\A\to P'$ into a partially
ordered set $P'$ satisfying the above conditions, we have
$$\m'(X)\leq\m'(Y)\quad\Longleftrightarrow\quad
\m(X)\leq\m(Y)\quad\text{for all}\quad X,Y\in \ind\A.$$
\end{thm}
\begin{proof} 
Use the axiomatic characterization of the chain length function
$\ell^*$ in Theorem~\ref{se:axiom}.
\end{proof}

\subsection*{Gabriel's main property}\label{se:main}
Let $\ell$ be a fixed length function on $\A$. The following main
property of the Gabriel-Roiter measure $\m=\ell^*$ is crucial; it is
the basis for all applications.

\begin{prop}[Gabriel]
Let $X,Y_1,\ldots,Y_r\in\ind\A$. Suppose that $X\subseteq
Y=\oplus_{i=1}^rY_i$. Then $\m(X)\leq\max\m(Y_i)$ and $X$ is a direct
summand of $Y$ if $\m(X)=\max\m(Y_i)$.
\end{prop}
\begin{proof}
The proof only uses the properties (GR1) -- (GR3) of $\m$.  Fix a
monomorphism $\p\colon X\to Y$. We proceed by induction on
$n=\ell_1(X)+\ell_1(Y)$.  If $n= 2$, then $\p$ is an isomorphism and the
assertion is clear. Now suppose $n>2$.  We can assume that
for each $i$ the $i$th component $\p_i\colon X\to Y_i$ of $\p$ is an
epimorphism. Otherwise choose for each $i$ a decomposition
$Y_i'=\oplus_j Y_{ij}$ of the image of $\p_i$ into indecomposables.
Then we use (GR1) and have $\m(X)\leq\max\m(Y_{ij})\leq\max\m(Y_i)$
because $\ell_1(X)+\ell_1(Y')<n$ and $Y_{ij}\subseteq Y_i$ for all $j$. Now
suppose that each $\p_i$ is an epimorphism. Thus
$\ell(X)\geq\ell(Y_i)$ for all $i$.  Let $X'\subset X$ be a proper
indecomposable subobject.  Then $\m(X')\leq\max\m(Y_i)$ because
$\ell_1(X')+\ell_1(Y)<n$, and $X'$ is a direct summand if
$\m(X')=\max\m(Y_i)$. We can exclude the case that
$\m(X')=\max\m(Y_i)$ because then $X'$ is a proper direct summand of
$X$, which is impossible. Now we apply (GR3) and obtain
$\m(X)\leq\max\m(Y_i)$.  Finally, suppose that
$\m(X)=\max\m(Y_i)=\m(Y_k)$ for some $k$. We claim that we can choose
$k$ such that $\p_k$ is an epimorphism.  Otherwise, replace all $Y_i$
with $\m(X)=\m(Y_i)$ by the image $Y_i'=\oplus_j Y_{ij}$ of $\p_i$ as
before.  We obtain $\m(X)\leq\max\m(Y_{ij})<\m(Y_k)$ since
$Y_{kj}\subset Y_k$ for all $j$, using (GR1) and (GR2).  This is a
contradiction.  Thus $\p_k$ is an epimorphism and in fact an
isomorphism because $\ell(X)=\ell(Y_k)$ by (GR2). In particular, $X$ is
a direct summand of $\oplus_iY_i$. This completes the proof.
\end{proof}

\subsection*{Gabriel-Roiter filtrations}\label{se:filt}
We keep a length function $\ell$ on $\A$ and the corresponding
Gabriel-Roiter measure $\m=\ell^*$.  Let $X,Y\in\ind\A$. We say that
$X$ is a {\em Gabriel-Roiter predecessor} of $Y$ if $X\subset Y$ and
$\m(X)=\max_{Y'\subset Y}\m(Y')$. Note that each object $Y\in\ind\A$ which is
not simple admits a Gabriel-Roiter predecessor because $\m$ is a
length function on $\ind\A$. A Gabriel-Roiter predecessor $X$ of $Y$
is usually not unique, but the value $\m(X)$ is determined by $\m(Y)$.

A sequence $$X_1\subset X_2\subset \ldots\subset X_{n-1} \subset X_n=X$$ in $\ind\A$ is called a
{\em Gabriel-Roiter filtration} of $X$ if $X_1$ is simple and
$X_{i-1}$ is a Gabriel-Roiter predecessor of $X_i$ for all $1< i\leq
n$. Clearly, each $X$ admits such a filtration and the values
$\m(X_i)$ are uniquely determined by $X$. Note that (C0) implies
\begin{equation}\label{eq:ell*}
\m(X)=\{\ell(X_i)\mid 1\leq i\leq n\}.
\end{equation}
\subsection*{Injective objects}
In order to illustrate Gabriel's main property, let us show that the
Gabriel-Roiter measure detects injective objects. This is a remarkable
fact because the Gabriel-Roiter measure is a combinatorial invariant,
depending only on the poset of indecomposable objects and some length
function, whereas the notion of injectivity involves all morphisms of
the category $\A$.

\begin{thm}\label{th:inj}
An indecomposable object $Q$ of $\A$ is injective if and only if there
is a length function $\ell$ on $\A$ such that $\ell^*(X)\leq\ell^*(Q)$
for all $X\in\ind\A$.
\end{thm}

We need the following lemma.

\begin{lem}\label{le:soc}
Let $\ell$ be a length function on $\A$ and fix indecomposable objects
$X,Y\in\A$. Suppose that for each pair of simple subobjects
$X'\subseteq X$ and $Y'\subseteq Y$, we have $\ell (X')<\ell (Y')$. Then 
$\ell^*(X)>\ell^*(Y)$.
\end{lem}
\begin{proof}
We choose Gabriel-Roiter filtrations $X_1\subset  \ldots \subset X_n=X$
and $Y_1\subset \ldots \subset Y_m=Y$. Then $\ell(X_1)<\ell(Y_1)$ and
the formula \eqref{eq:ell*} implies $$\ell^*(X)=\{\ell(X_i)\mid 1\leq
i\leq n\}>\{\ell(Y_i)\mid 1\leq i\leq m\}=\ell^*(Y).$$
\end{proof}

\begin{proof}[Proof of the theorem]
Suppose first that $Q$ is injective. Then $Q$ has a unique simple
subobject $S$ and we define a length function $\ell=\ell_S$ on $\A$ by
specifying its values on each simple object $T\in\A$ as follows:
$$\ell(T):=\begin{cases} 1 \quad \text{if } T\cong S,\\ 2 \quad
\text{if } T\not\cong S.\end{cases}$$ Now let $X\in\ind\A$. We claim
that $\ell^*(X)\leq\ell^*(Q)$. To see this, let $X'\subseteq X$ be the
maximal subobject of $X$ having composition factors isomorphic to
$S$. Using induction on the composition length $n=\ell_1(X')$ of $X'$,
one obtains a monomorphism $X'\to Q^n$, and this extends to a map
$\p\colon X\to Q^n$, since $Q$ is injective. Let $X/X'=\oplus_iY_i$ be
a decomposition into indecomposables and $\pi\colon X\to X/X'$ be the
canonical map. Note that $\ell^*(Y_i)<\ell^*(Q)$ for all $i$ by our
construction and Lemma~\ref{le:soc}.  Then $(\pi,\p)\colon X\to
(\oplus_iY_i) \oplus Q^n$ is a monomorphism and therefore
$\ell^*(X)\leq \ell^*(Q)$ by the main property.

Suppose now that $\ell^*(X)\leq\ell^*(Q)$ for all $X\in\ind\A$ and
some length function $\ell$ on $\A$. To show that $Q$ is injective,
suppose that $Q\subseteq Y$ is the subobject of some $Y\in\A$. Let
$Y=\oplus Y_i$ be a decomposition into indecomposables. Then the main
property implies $\ell^*(Q)\leq\max\ell^*(Y_i)\leq\ell^*(Q)$ and
therefore $Q$ is a direct summand of $Y$. Thus $Q$ is injective and
the proof is complete.
\end{proof}

Let us mention that there is the following analogous characterization
of the simple objects of $\A$.

\begin{cor}
An indecomposable object $S$ of $\A$ is simple if and only if there is
a length function $\ell$ on $\A$ such that $\ell^*(S)\leq\ell^*(X)$
for all $X\in\ind\A$.
\end{cor}
\begin{proof}
Use the property (GR1) of the Gabriel-Roiter measure and apply
Lemma~\ref{le:soc}.
\end{proof}

\subsection*{The Kronecker algebra}
Let $\La=\smatrix{k&k^2\\ 0&k}$ be the Kronecker algebra over an
algebraically closed field $k$.  We consider the abelian length
category which is formed by all finite dimensional $\La$-modules. A
complete list of indecomposable objects is given by the preprojectives
$P_n$, the regulars $R_n(\a,\b)$, and the preinjectives $Q_n$. More
precisely,
$$\ind\La=\{P_n\mid n\in\bbN\}\cup\{R_n(\a,\b)\mid
n\in\bbN,\,(\a,\b)\in\mathbb P_k^1\}\cup\{Q_n\mid n\in\bbN\},$$ and we
obtain the following Hasse diagram.
$$
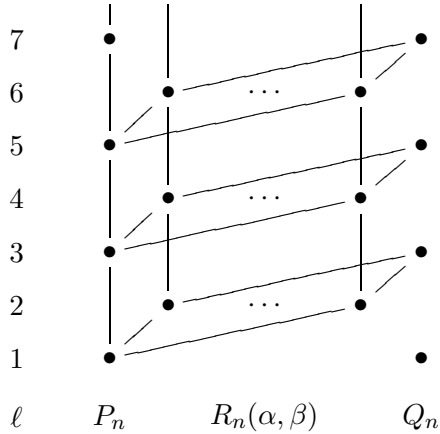
\xymatrix@=0.6em{
&&\ar@{-}[d]&\ar@{-}[dd]&&\ar@{-}[dd]\\
7&&\bullet\ar@{-}[dd]&&&&\bullet\ar@{-}[dl]\ar@{-}[dlll]\\
6&&&\bullet\ar@{-}[dd]\ar@{-}[dl]&\cdots&\bullet\ar@{-}[dd]\ar@{-}[dlll]\\
5&&\bullet\ar@{-}[dd]&&&&\bullet\ar@{-}[dl]\ar@{-}[dlll]\\
4&&&\bullet\ar@{-}[dd]\ar@{-}[dl]&\cdots&\bullet\ar@{-}[dd]\ar@{-}[dlll]\\
3&&\bullet\ar@{-}[dd]&&&&\bullet\ar@{-}[dl]\ar@{-}[dlll]\\
2&&&\bullet\ar@{-}[dl]&\cdots&\bullet\ar@{-}[dlll]\\
1&&\bullet&&&&\bullet\\ \ell&&P_n&&R_n(\a,\b)&&Q_n }$$ The set of
indecomposables is ordered as follows via the Gabriel-Roiter measure
with respect to $\ell=\ell_1$.
\begin{align*}
\ell^*&:\quad Q_1=P_1<P_2<P_3<\ldots \;\;\;R_1<R_2<R_3<\ldots
\;\;\;\ldots<Q_4<Q_3<Q_2\\
(\ell^*)^*&:\quad Q_1=P_1<R_1<Q_2<P_2<R_2<Q_3<P_3<R_3<Q_4<\ldots 
\end{align*}
Moreover, $((\ell^*)^*)^*$ and $\ell^*$ are equivalent length
functions.

\begin{rem}
While $\ell^*$ has been successfully employed for
proving the first Brauer-Thrall conjecture, Hubery points out that
$(\ell^*)^*$ might be useful for proving the second.  In fact, one
needs to find a value $(\ell^*)^*(X)$ such that the set
$\{X'\in\ind\La\mid (\ell^*)^*(X')=(\ell^*)^*(X)\}$ is infinite. The
example of the Kronecker algebra shows that there exists such a value
having only finitely many predecessors
$(\ell^*)^*(Y)<(\ell^*)^*(X)$. Note that in all known examples
$((\ell^*)^*)^*$ and $\ell^*$ are equivalent.
\end{rem}

\subsection*{Acknowledgements} 
This material has been presented at the ``Advanced School and
Conference on Representation Theory and Related Topics'' in Trieste
(ICTP, January 2006) and I am grateful to the organizers. In addition,
I wish to thank Philipp Fahr, Andrew Hubery, Osamu Iyama, and Karsten
Schmidt for helpful discussions and comments.


\begin{thebibliography}{99}
%
\bibitem{G}{\sc P.\ Gabriel:} Indecomposable representations
II. Symposia Mathematica {\bf 11} (1973), 81--104.
%
\bibitem{R1}{\sc C.\ M.\ Ringel:} The Gabriel-Roiter
measure.  Bull. Sci. Math. {\bf 129}  (2005), 726--748.
%
\bibitem{R2}{\sc C.\ M.\ Ringel:} Foundation of the representation
theory of Artin algebras, using the Gabriel-Roiter measure. In: Trends
in Representation Theory of Algebras and Related
Topics. Contemp. Math. {\bf 406} (2006), 105--135.
%
\bibitem{R3}{\sc C.\ M.\ Ringel:} The theorem of Bo Chen and Hall
polynomials. Nagoya Math. J. {\bf 183} (2006).
%
\bibitem{Ro}{\sc A.\ V.\ Roiter:} Unboundedness of the dimension of the
indecomposable representations of an algebra which has infinitely many
indecomposable representations. Izv. Akad. Nauk SSSR Ser. Mat. {\bf
32} (1968), 1275-1282.
%
\end{thebibliography}
\end{document}